\documentstyle[12pt]{article}
\title{A stratification of the moduli space of vector bundles on curves}
\author{L.~Brambila-Paz 
\and
Herbert Lange\thanks{The first author is member of the VBAC 
Research group of Europroj. The work forms part of the project 
``Moduli of bundles in algebraic geometry", funded by 
the EU International Scientific Cooperation Initiative (Contract
No.~CI1*-CT93-0031). The first author also acknowledges 
support from CONACYT (Project No. 3231-E9307).}}
\date{ }

\def\z{{\bf Z}}
\def\p{{\bf P}}
\def\se{{\it s}}
\def\ker{\mbox{\rm ker}}

\def\mod{\mbox{\rm mod}}
\let\dps=\displaystyle
\def\m@th{\mathsurround=\z@}

\newcommand{\qed}{\hspace{\fill}$\Box$}

\newtheorem{thm}{Theorem}[section]

\newtheorem{propn}[thm]{Proposition}
\newtheorem{rem}[thm]{Remark}
\newtheorem{lemma}[thm]{Lemma}
\newtheorem{cor}[thm]{Corollary}

\newcommand{\ra}{\rightarrow}

\begin{document}
\maketitle

\rightline{{\it Dedicated to M. Kneser}}

\bigskip

\bigskip

{\bf Introduction}

Let $E$ be a vector bundle of rank 2 on a smooth projective curve $C$ of
genus $g \geq 2$ over an algebraically closed field $K$ of arbitrary
characteristic.

The invariant $${\se}_1(E):=\deg E-2 \max\deg(L),$$ where the maximum is
taken over all line subbundles $L$ of $E$, is just the minimum of the self
intersection numbers of all sections of the ruled surface $\p(E)\to C$.
Note that $E$ is stable (respectively semistable) if and only if
${\se}_1(E)\ge 1$ (respectively $\ge 0$). According to a Theorem of C.
Segre ${\se}_1(E)\le g$. Moreover, the function ${\se}_1$ is lower
semicontinuous. Thus ${\se}_1$ gives a stratification of the moduli space
${\cal M}(2,d)$ of stable vector bundles of rank $2$ and degree $d$  on
$C$, into locally closed subsets ${\cal M}(2,d,s)$ according to the value
$s$ of ${\se}_1(E)$.

It is shown in \cite{ln} that for $s>0$, $s\equiv d\; \mod\; 2$ the
algebraic variety ${\cal M}(2,d,s)$ is non-empty, irreducible and of
dimension $$\dim {\cal M}(2,d,s)=\left\{\begin{array}{lcl} 3g+s-2 &  &
s\le g-2\\ &{\rm if}& \\ 4g-3 & & s\ge g-1\end{array}\right. \eqno(A)$$
Let ${M}_1(E)$ denote the set of line subbundles of $E$ of maximal degree.
The set  ${M}_1(E)$ can be considered as an algebraic scheme in a natural
way. Maruyama proved in \cite{ma} (see also \cite{ln} for different
proofs) the following statement:

For a general vector bundle $E$ in ${\cal M}(2,d,s)$: $$\dim
{M}_1(E)=\left\{\begin{array}{lcl} 1 &  & s=g\\ &{\rm if}& \\ 0& & s\le
g-1\end{array}\right.\eqno(B)$$

It is the aim of the present paper, to generalize statements (A) and (B)
to vector bundles of arbitrary rank $r\ge 2$.

Let now $E$ be a vector bundle of rank $r\ge 2$ over $C$. For any integer
with $1\le k\le r-1$ we define $${\se}_k(E):=k\deg E-r\max\deg F.$$ where
the maximum is taken over all subbundles $F$ of rank $k$ of $E$.  There is
also an interpretation of the invariant ${\se}_k(E)$ in terms of self
intersection numbers in the ruled variety $\p(E)$ (see \cite{l3}).
According to a theorem of Hirschowitz (see \cite{h1}) $$\se_k(E)\le
k(r-k)(g-1)+(r-1).$$

Moreover, the function ${\se}_k$ is lower semicontinuous. Thus ${s}_k$
gives a stratification of the moduli space ${\cal M}(r,d)$ of stable
vector bundles of rank $r$ and degree on $d$ on $C$ into locally closed
subsets ${\cal M}(r,d,k,s)$ according to the value of $s$ and $k$. There
is a component ${\cal M}^0(r,d,k,s)$ of ${\cal M}(r,d,k,s)$  distinguish
by the fact that a general $E\in {\cal M}^0(r,d,k,s)$ admits a stable
subbundle $F$ such that $E/F$ is also stable.

In this paper we prove (see Theorem 4.2):

{\it For $g\ge \frac{r+1}{2}$ and $0<s\leq k(r-k)(g-1) +(r+1)$, $s\equiv
kd\;\mod\; r,$ ${\cal M}^0(r,d,k,s)$ is non-empty, and its component
${\cal M}^0(r,d,k,s)$ is of dimension } $$\dim {\cal
M}^0(r,d,k,s)=\left\{\begin{array}{lcl} (r^2+k^2-rk)(g-1)+s-1&
&s<k(r-k)(g-1)\\ &{\rm if}&\\ r^2(g-1)+1& & s\ge
k(r-k)(g-1)\end{array}\right.\eqno(C)$$

The bound $g \geq \frac{r+1}{2}$ works for all $k, 1 \leq k \leq r-1$
simultaneously. For some special $k$ the result is better. For example,
for $k = 1$ or $r-1$ statement (C) is valid for all $g \geq 2$ (see Remark
3.3).

Let $M_k(E)$ denote the set of maximal subbundles of rank $k$ of $E$  and
${\widehat{M_k}}(E)$ denote the closure of the set of stable subbundles
$F$ of rank $k$ of $E$ of maximal degree such that $E/F$ is also stable in
$M_k(E)$. Also ${M}_k(E)$ can be considered as an algebraic scheme in a
natural way and ${\widehat{M_k}}(E)$ is a union of components of $M_k(E).$
 
Theorem 4.4 below says:

{\it For a general vector bundle $E$ in ${\cal M}^0(r,d,k,s)$ we have }
$$\dim {M}_k(E)=\max(s-k(r-k)(g-1),0)\eqno(D)$$ {\it again under the
hypothesis} $g\ge \frac{r+1}{2}$.

Note that (C) and (D) are direct generalizations of (A) and (B).

The main difficulty in the proofs is to show that ${\cal M}(r,d,k,s)$ is
nonempty.  (see Theorem 3.2 below). This was shown in \cite{l1} and
\cite{h2} for $s\ge k(r-k)(g-1)$ i.e. in the generic case and for
$s\le{k(r-k)(g-1)\over \max(k,r-k)}$ in
\cite{bbr}. For a generic curve of genus $g \geq 2$ this was proven by M.
Teixidor in [16].  Special cases were also considered by B. Russo in
\cite{br}, (see also \cite{bb}) and M. Teixidor in \cite{tei}.  The idea
of proof is to start with a general vector bundle $E_0$ with
${\se}_k(E_0)\ge k(r-k)(g-1)$, then construct a sequence of elementary
transformations to find a vector bundle in $M(r,d,k,s)$.  Statement $(D)$
is an easy consequence of statement $(C)$.
\bigskip

{\bf Acknowledgement} Part of the work was done during a visit of both
authors to CIMAT, Gto. Mexico. They wish to acknowledge the hospitality of
CIMAT. We would like to thank P.E. Newstead for his comments on the first
version of this paper and Remark 4.5.  The first author also want to thank
P. E. Newstead and B. Russo for previous discussions on the subject and
the University of Liverpool where those discussion were done.

\section{The invariants ${ {}{\se_k}}(E)$}

Let $C$ be a smooth projective curve of genus $g\ge 2$ over an
algebraically closed field $K$ of arbitrary characteristic.  and let $E$
denote a vector bundle of rank $r\ge 2$ over $C$. For any integer $k$ with
$1\le k\le r-1$ let ${}{Sb_k}(E)$ denote the {\it set of subbundles of
rank $k$ of $E$}.  If we denote by $\xi$ the generic point of the curve
$C$, then it is easy to see that there is a canonical bijection between
${}{Sb_k}(E)$ and the set of $k$-dimensional subvector spaces of the
$K(\xi)$-vector space $E(\xi)$.  For any subbundle $F\in {}{Sb_k}(E)$
define the integer ${\se_k}(E,F)$ by $${\se_k}(E,F):= k\deg E-r\deg \;F.$$
The vector bundle $E$ does not admit subbundles of arbitrarily high
degree. Hence $${}{\se_k}(E):=\mathop{\min}_{F\in
{}{Sb_k}(E)}{\se_k}(E,F)$$ is a well defined integer depending only  on
$E$ and $k$.

\begin{rem} \label{1.1} \begin{em}
 The slope of a vector bundle $F$ on $C$ is defined as $\mu(F)={\deg\;
F\over {\rm rk} F}$. If $F$ is a subbundle of rank $k$ of $E$ then
$${\se_k}(E,F)=k(r-k)\left(\mu(E/F)-\mu(F)\right).$$
\end{em} \end{rem}
In particular $${}{\se_k}(E)=k(r-k)\cdot \mathop{\min}_{F\in
{}{Sb_k}(E)}\left(\mu(E/F)-\mu(F)\right)$$ So instead of the invariant
${}{\se_k}(E)$ one could also work with the invariant \par
\noindent{$\min_{F\in {}{Sb_k(E)}}\left(\mu(E/F)-\mu(F)\right).$} 
However, for some proofs it is more convenient to work with integers. Note
that there is also a geometric interpretation of the invariant
${}{\se_k}(E)$ in terms of intersection numbers on the associated
projective bundle $\p(E)$ (see \cite{l3}).

\begin{rem} \label{1.2} \begin{em} The following properties of the
invariant ${}{\se_k}(E)$ are easy to see (see \cite{l1})\begin{itemize}
\item[(a)] ${}{\se_k}(E\otimes L)={}{\se_k}(E)$ for all $L\in Pic(C)$.
\item[(b)] ${}{\se_k}(E)={}{\se_{r-k}}(E^*)$.
\item[(c)] $E$ is stable (respectively semistable) if and only if
${}{\se_k}(E)>0$ 
(respectively ${}{\se_k}(E)\ge 0$) for all $1\le k\le r-1$.
\item[(d)] Let $T$ be an algebraic scheme over $K$ and ${\cal E}$ a vector
bundle of rank $r$ on $C\times T$. For any point $t\in T,$ let $\bar{t}$
denote a geometric point over $t$.  The function ${}{s_k}:T\to\z $ defined
as $t\mapsto{}{\se}_k({\cal E} |_{C\times\{\bar{t}\}})$ is well defined
and lower semicontinuous.
\end{itemize}
\end{em} \end{rem}

Whereas the function ${}{s_k}$ may take arbitrarily negative values (for
suitable direct sums of line bundles).  However, it is shown in \cite{ms}
and \cite{l1} that ${}{\se_k}(E)\le k(r-k)g$. Hirschowitz gives in
\cite{h1} the better bound, $${}{\se_k}(E)\le k(r-k)(g-1)+(r-1).$$

We want to study the behaviour of the invariant ${}{\se_k}(E)$ under an
elementary transformation of the vector bundle $E$. Recall that an {\it
elementary transformation} $E'$ of $E$ is defined by an exact sequence
$$0\to E'\to E\mathop{\to}^{\ell}K(x)\to 0\eqno(1)$$ where $K(x)$ denotes
the skyscraper sheaf with support $x\in C$ and fibre $K$. Since $\ell$
factorizes uniquenly via a $k$-linear form $E(x)\stackrel{\ell}\rightarrow
K(x)$, also denoted by $\ell$, the set of elementary transformations of
$E$ is parametrized by pairs $(x,\ell)$ where $x$ is a closed  point of
$C$ and $\ell$ is a linear form on the vector space $E(x)$.

So the set of elementary transformations of $E$, which we denote by
$elm(E)$, forms a vector bundle of rank $r$ over the curve $C$. Note that
for any $E'\in\;  elm(E)$ $${\rm rk} (E')={\rm rk}(E)\quad{\rm
and}\quad\deg E'=\deg E-1.$$

\begin{lemma} \label{1.3}
 For any $E'\in elm(E)$ the map $\varphi:{}{Sb_k}(E)\to {}{Sb_k}(E') $
defined by $F\mapsto F\cap E'$ is a bijection.
\end{lemma}
\noindent{\bf Proof.} For the proof only note that the inverse map is given as 
follows: Suppose $F'\in {}{Sb_k}(E')$. Consider $F'$ as a subsheaf of $E$
and let $F$ denote the subbundle of $E$ generated by $F'$.  The map
$F'\mapsto F$ is inverse to $\varphi$.\qed

Now consider a subbundle $F\in {}{Sb_k}(E)$ and denote $F'=\varphi(F)\in
{}{Sb_k}(E')$.  In order to compute the number ${\se_k}(E',F')$ we have to
distinguish two cases.  We say that the subbundle $F$ is of {\it type I}
with respect to $E'$ if $F\subseteq E'$ and $F$ is of {\it type II} with
respect to $E'$ otherwise. Let $(x,\ell)\in elm(E)$ denote the pair
defining the elementary transformation $E'$ of $E$. We obviously have

\begin{lemma} \label{1.4}
The subbundle $F\in {}{Sb_k}(E)$ is of type I with respect to $E'$ if and
only if the linear form $\ell :E(x)\to K(x)$ vanishes on the subvector
space $F(x)$ of $E(x)$.
\end{lemma}

\begin{lemma} \label{1.5} 
If $E'$ is an elementary transformation of $E$, $F\in {}{Sb_k}(E)$ and
$F'=\varphi(F)$ then, \begin{itemize}
\item[(i)] ${\se}_k(E',F') ={\se}_k(E,F)-k$   if $F$ is of type I with
respect to $E'$
\item[(ii)] ${\se}_k(E',F')={\se}_k(E,F)+(r-k)$   if $F$ is of type II
with respect to $E'$.
\end{itemize}
\end{lemma}
\noindent{\bf Proof.} If $F$ is of type I, we have the following diagram
$$\begin{array}{ccccccccc} &     &   &    &  0         &    &  0        &
&    \\ &     &   &    &\downarrow  &    &\downarrow &    &    \\ 0 & \to
& F'&\to & E'         &\to & E'/F'     &\to & 0 \\ &     & \|&
&\downarrow  &    &\downarrow &    &    \\ 0 & \to & F &\to & E
&\to & E/F     &\to & 0 \\ &     &   &    &\downarrow  &    &\downarrow &
&    \\ &     &   &    &K(x) & =    &K(x)  &    &    \\ &     &   &
&\downarrow  &    &\downarrow &    &    \\ &     &   &    &  0         &
&  0        &    &    \\
\end{array}$$
Hence ${\se}_k(E',F')=k(\deg E-1)-r\deg F={\se}_k(E,F)-k$.  If $F$ is type
II we have the following diagram $$\begin{array}{ccccccccc} &     &   0
&    &  0         &        &           &    &    \\ &     & \downarrow &
&\downarrow  &        &           &    &    \\ 0 & \to & F'         &\to &
E'         &\to     & E'/F'     &\to & 0 \\ &     &\downarrow  &
&\downarrow  &        &\| &&    \\ 0 & \to & F          &\to & E
&\to     & E/F       &\to & 0 \\ &     & \downarrow &    &\downarrow  &
&           &    &    \\ &     & K(x)       & =  & K(x)       &        &
&    &    \\ &     &\downarrow  &    &\downarrow  &        &           &
&    \\ &     &      0     &    &  0         &        &           &    &
\\
\end{array}$$
Hence ${\se}_k(E',F')=k(\deg E-1)-r(\deg F-1)= {\se}_k(E,F)+(r-k).$\qed

A {\it maximal subbundle} $F\in {}{Sb_k}(E)$ is by definition a  subbundle
of rank $k$ of maximal degree of $E$. Note that $F\in {}{Sb_k}(E)$ is a
maximal subbundle if and only if $${}{\se_k}(E)={\se}_k(E,F).$$ An
elementary transformation $E'$ of $E$ will be called {\it of $k$-type I}
if $E$ admits a maximal subbundle of rank $k$ which is of type I with
respect to $E'$. Otherwise $E'$ will be called {\it of $k$-type II}.

\begin{propn} \label{1.6} If $E'$ is an elementary transformation of $E$, then 
$${}{\se_k}(E')=\left\{\begin{array}{lcl} {}{\se_k}(E)-k& &E'\quad\mbox{is
of $k$-type I} \\ &{\rm if} & \\ {}{\se_k}(E)+(r-k)& &E'\quad \mbox{is of
$k$-type II. }\end{array}\right. $$
\end{propn}
\noindent{\bf Proof.} Let $F\in {}{Sb_k}(E)$ and 
$F'=\varphi(F)\in {}{Sb_k}(E')$. Suppose first $E'$ is of $k$-type I. If
$F$ is maximal and of type I with respect to $E'$, then
${\se_k}(E',F')={}{\se_k}(E)-k$.  If $F$ is maximal and of type II, then
${\se_k}(E',F')={}{\se_k}(E)+(r-k)$.  If $F$ is not maximal, then $\deg \;
F\le{1\over r}(k\deg E-{}{\se_k}(E))-1$ and so $${\se}_k(E,F)\ge k\deg E-
r\deg F\ge{}{\se_k}(E)+r\eqno(2)$$ Hence
${\se}_k(E',F')\ge{\se}_k(E,F)-k\ge{}{\se_k}(E)+(r-k)$. This implies the
assertion if $E'$ is of $k$-type I.

Suppose now $E'$ is of $k$-type II. Any maximal subbundle $F \subset E$ is
by assumption of type II with respect to $E'$. Hence according to Lemma
1.5 ${\se}_k(E',F')={}{\se_k}(E)+(r-k)$. If $F$ is not maximal, then (2)
and Lemma 1.5 imply
$${\se}_k(E',F')\ge{\se}_k(E,F)-k\ge{}{\se_k}(E)+(r-k).$$\qed

\begin{rem} \label{1.7} \begin{em}
 One inmediately deduces from the proof of Proposition 1.6:\begin{itemize}
\item[(i)] If $E'$ is 
of $k$-type I, then the maximal subbundles of rank $k$ of $E'$ are exactly
the maximal subbundles of rank $k$ of $E$ which are of type I with respect
to $E'$.
\item[(ii)] If $E'$ is of $k$-type II, then the maximal subbundles 
of rank $k$ of $E'$ are exactly the subbundles $F'=\varphi(F)$, where
$F\in {}{Sb_k}(E)$ is either maximal or of degree one less than the degree
of a maximal subbundle and of type I with respect to $E'$.
\end{itemize}
\end{em}\end{rem}
Dualizing the exact sequence (1) we obtain an exact sequence $$0\to E^*\to
E'^*\to K(x)\to 0$$ Hence $E^*$ is an elementary transformation of $E'^*$,
called the {\it dual elementary transformation}.

\begin{cor} \label{1.8}
 For an elementary transformation $E'$ of $E$ the following conditions are
equivalent\begin{itemize}
\item[(i)] $E'$ is of $k$-type I.
\item[(ii)] The dual elementary transformation $E^*$ of $E'^*$ is of
$(r-k)$-type II.
\end{itemize}
\end{cor}

\noindent{\bf Proof:} According to Proposition 1.6 and Remark 1.7, (i)
holds if 
and only if 
${}{\se_k}(E')={}{\se_k}(E)-k$. But 
${}{\se_{r-k}}(E^*)={}{\se_k}(E)$ (see Remark 1.2, (b)). 
Hence (i) holds if and only if 
${}{\se_{r-k}}(E'^*)={}{\se_{r-k}}(E^*)-k$ i.e. if and only if
${}{\se_{r-k}}(E^*) ={}{\se_{r-k}}(E'^*)+ r-(r-k)$. Applying 
Proposition 1.6 again 
gives the assertion.\qed

\section{Maximal subbundles}

Let $E$ denote a vector bundle of rank $r$ and degree $d$ on the curve
$C$.  In this section we study the set ${}{M_k}(E)$ of maximal subbundles
of rank $k$ of $E$. Let $d_k$ denote the common degree of the maximal
subbundles of rank $k$ of $E$. The following lemma shows that ${}{M_k}(E)$
admits a natural structure of a projective scheme over $K$. Denote by
$Q:={\rm Quot}^{r-k, d-d_k}_E$ the Quot scheme of coherent quotients of
rank $r-k$ and degree $d-d_k$ of $E$.
\begin{lemma} \label{2.1} 
 There is a canonical identification of ${}{M_k}(E)$ with the set of
closed points of $Q$.
\end{lemma}
\noindent{\bf Proof.} If $F\in {}{M_k}(E)$, then $E\to E/F$ gives a closed
point of $Q$.
On the other hand if $\dps E\mathop{\to}^p G\to 0$ corresponds to a closed
point of $Q$, then $F=\ker \; p\in {}{M_k}(E).$\qed

Let ${\cal G}$ denote the universal quotient sheaf on $C\times Q$.  The
maximality condition implies that ${\cal G}$ is locally free.  Hence if
$\dps Grass_{r-k}(E)\mathop{\to}^p C$ denotes the Grassmanian scheme of
$(r-k)$-dimensional quotient vector spaces of the fibres $E(x)$ and
$p^*E\to{\cal U}\to 0$ the universal quotient on $Grass_{r-k}(E)$, then
any $F\in {}{M_k}(E)$ corresponds on the one hand to a closed point $t$ of
$Q$ and on the other hand to a section $\sigma_t: C\to Grass_{r-k}(E)$.
This leads to a morphism $$\phi :\left\{\begin{array}{l} C\times Q\to
Grass_{r-k}(E)\\ (x,t)\mapsto \sigma_t(x)\end{array}\right. $$ with the
property that ${\cal G}=\phi^*{\cal U}$.

\begin{lemma} \label{2.2} The morphism $\phi :C\times Q\to Grass_{r-k}(E)$
is finite.
\end{lemma}
For a proof we refers to \cite{ms} or \cite{l3}, Lemma 3.9.

The geometric interpretation of Lemma 2.2 is

\begin{cor} \label{2.3} 
 Let $x\in C$ and $V\subset E(x)$ and $k$-dimensional subvector space.
There are most finitely many maximal subbundles $F$ of rank $k$ of $E$
such that $F(x)=V$.
\end{cor}

\begin{cor} \label{2.4} dim ${}{M_k}(E)\le k(r-k).$
\end{cor}
\noindent{\bf Proof:} From Lemma 2.1 there is a canonical 
identification of ${}{M_k}(E)$ with the set of closed points of $Q$.
Since $\phi$ is finite, we have that dim $Q\le {\rm dim}\
Grass_{r-k}(E)-$dim$ C=k(r-k)+1-1.$\qed

Assume now that dim ${}{M_k}(E)=n$, where $n\le k(r-k)$ according to
Corollary 2.4, and let $E'$ be an elementary transformation of $E$. We
want to estimate dim ${}{M_k}(E')$.  Suppose that $E'$ corresponds to the
pair $(x,\ell)$ with exact sequence (1) of Section 1. For any
$(r-1)$-dimensional subvector space $V$ of the vector space $E(x)$
consider the Schubert cycle $$\sigma_k(V):=\{F\in Grass_{r-k}(E(x))|
F\subset V\}.$$
\begin{propn} \label{2.5}If $E'$ is of $k$-type I, then 
dim ${}{M_k}(E')\ge {\rm dim} {}{M_k}(E)-k$.
\end{propn}
\noindent{\bf Proof.} Denote $V=\ker(\ell :E(x)\to K)$. According to 
Remark 1.7, there is a canonical identification $${}{M_k}(E')=\{F\in
{}{M_k}(E)|F(x)\subseteq V\}.$$ Defining $${}{M_k}(E)(x):=\{F(x)|F\in
{}{M_k}(E)\},$$ we have from Corollary 2.3
\begin{eqnarray*}
\dim {}{M_k}(E')&=&\dim {}{M_k}(E')(x)\\
&=&\dim ({}{M_k}(E)(x)\cap\sigma_k(V))\\ &\ge&\dim
{}{M_k}(E)(x)+\dim\sigma_k(V)-\dim Grass_{r-k}(E(x))\\ &=&
n+k(r-1-k)-k(r-k)=n-k. \hspace{5cm} \Box    \end{eqnarray*} This completes
the proof of the proposition.\qed

\begin{propn} \label{2.6} 
If $E'\in elm(E)$ is of $k$-type II, then {\rm dim} ${}{M_k}(E')\le ${\rm
dim} ${}{M_k}(E)+r-k$.
\end{propn}

\noindent{\bf Proof.} According to Corollary 1.8, the dual elementary 
transformation $E^*$ of $E'^*$ is of $(r-k)$-type I. Hence by Proposition
2.5 $$\dim {}{M_{r-k}}(E^*)\ge\dim {}{M_{r-k}}(E'^*)-(r-k).$$ But
dualizing induces a canonical isomorphism $\dps
{}{M_k}(E)\mathop{\to}^{\sim} {}{M_{r-k}}(E^*)$ and similarly for $E'$.
This completes the assertion.\qed

\section{Stable extensions}

Suppose now that the genus $g$ of the curve $C$ is $\ge \frac{r+1}{2}$.
Let $r,d,k$ and $s$ be integers with $r\ge 2$, $1\le k\le r-1$, $0<s\leq
k(r-k)(g-1)+(r+1)$ and $s\equiv kd\; \mod\; r$. The aim of this section is
the proof of the following theorem.

\begin{thm} \label{3.1} 
There exists an extension $0\to F\to E\to G\to 0$ of vector bundles on the
curve $C$ with the following properties \begin{itemize}
\item[(i)] rk $E=r$, deg $E=d$.
\item[(ii)] $F$ is a maximal subbundle of rank $k$ of $E$ with
${\se}_k(E,F)=s$.
\item[(iii)] $E,F$ and $G$ are stable.
\end{itemize}
\end{thm}

Let $d_1$ be the unique integer with $s=kd -rd_1$ and $d_2=d-d_1$.
According to
\cite{nr} Proposition 2.4 there are finite \'etale coverings
$$\pi_1:{\widetilde{M_1}}\to {\cal M}(k,d_1) \ \ {\rm and} \ \
\pi_2:\widetilde{M_2}\to {\cal M}(r-k, d_2)$$
such that there are vector bundles ${\cal F}_i$ on
$C\times\widetilde{M}_i$ whose classifying map is just $id \times \pi_i$
for $i=1,2$.  Let $p_{ij}$ denote the canonical projections of
$C\times\widetilde{M}_1\times \widetilde{M}_2$ for $i,j=0,1,2$. According
to \cite{l2}, Lemma 4.1 the sheaf $R^1p_{12*}(p^*_{02}{\cal F}^*_2\otimes
p_{01}^*{\cal F}_1)$ is locally free of rank $k(r-k)(g-1)+s$ on
$\widetilde{M}_1\times\widetilde{M}_2$. Let $$\pi
:\p:=\p(R^1p_{12*}(p^*_{02}{\cal F}_2^*\otimes p_{01}^*{\cal
F}_1)^*)\to\widetilde{M}_1\times\widetilde{M}_2$$ denote the corresponding
projective bundle. According to \cite{l2}, Corollary 4.5 there is an exact
sequence $$0\to\pi^*p^*_{01}{\cal F}_1\otimes {\cal O}_{\p}(1)\to{\cal
E}\to\pi^*p^*_{02} {\cal F}_2\to 0\eqno(3)$$ on $C\times\p $, universal in
a sense which is outlined in that paper.  In particular this means that
for every closed point $q\in \p$ the restriction of the exact sequence (3)
to $C\times\{q\}$ is just the extension of ${\cal
F}_2|_{C\times\{p_2(q)\}}$ by ${\cal F}_1|_{C\times\{p_1(q)\}}$ modulo
$K^*$,  which is represented by the point $q$. Here $p_i: \p\to
\widetilde{M}_i$ denotes the canonical map.

With $r,k,d$ and $s$ as above consider the set
$$U(r,d,k,s):=\{q\in\p:{\cal E} |_{C\times\{q\}}\;\;\mbox{is stable
with}\;\;{}{\se}_k({\cal E}|_{C\times\{q\}})=s\}$$

 From the lower semicontinuity of the function ${}{\se}_k$ and stability
being an open conditon we deduce that the set $U(r,d,k,s)$ is an open
subset of $\p$. Hence, Theorem 3.1 is equivalent to the following theorem.

\bigskip
\begin{thm} \label{3.2}
 For any $r,k,d$ and $s$ as above the set $U(r,d,k,s)$ is nonempty.
\end{thm}

\noindent{\bf Proof.} It suffices to show that the set 
$$U(r,d,k,s,i) :=\{q\in \p :s_i({\cal E} _{|_{C\times \{q\}}}) >0 \ \
\mbox{ and} \ \ s_k({\cal E} _{|_{C\times \{q\}}})=s \}$$ is nonempty for
any $i=1,...,r-1, i\not= k$, since $$U(r,d,k,s) =
\bigcap_{\stackrel{i=1}{i\not=k}}^{r-1} U(r,d,k,s,i) $$ and the function
$s_i$ is lower semicontinuous. According to Remark 1.2 (b) dualization
gives a canonical bijection $$U(r,d,k,s,i) \stackrel{\sim}{\ra}
U(r,-d,r-k,s,r-i).$$ Hence it suffices to show that $$ U(r,d,k,s,i) \not=
\emptyset,$$ for all $r,d,k,s $ as above and all $1\leq i\leq k-1.$ Choose
a positive integer $N_k$ such that $$k(r-k)(g-1) \leq s +N_kk\leq
k(r-k)(g-1) +r-1 \eqno(4)$$ and denote $\widetilde{d} :=d+N_k.$

We call a vector bundle $E$ out of the moduli space ${\cal
M}(r,\widetilde{d} )$ {\it general}, if for all $0<j< r$ the number
${}{\se_j}(E)$ takes a maximal value, say ${\se_{j,\max}}$. By the
semicontinuity of the function ${}{\se_j}$ the set of general vector
bundles is open and dense in $M(r,\widetilde{d})$.  According to a theorem
of Hirschowitz (see \cite{h1}, Th\'eor\`eme p. 153):
$${\se}_{j,\max}=j(r-j)(g-1)+\epsilon_j$$ where $\epsilon_j$ is the unique
integer with $0\le\epsilon_j\le r-1$ such that
$j(r-j)(g-1)+\epsilon_j\equiv j{\widetilde{d}}\;\mod\; r$. Moreover, it is
shown in \cite{l1} (p. 458), that $U(r,\widetilde{d},k,s_{k,\max})$ is
non-empty and its image is  open and dense in $M(r,\widetilde{d})$ for all
$r,\widetilde{d}$ and $k$.

Let $0\ra F_0 \ra E_0 \ra G_0 \ra 0 $ be an exact sequence corresponding
to a general point in $U(r,\widetilde{d},k,s_{k,\max}).$ Then $E_0, F_0$
and $G_0$ respectively are general vector bundles in ${\cal
M}(r,\widetilde{d}), {\cal M}(k,d_k)$ and ${\cal M}(r-k,\widetilde{d}
-d_k)$ respectively, with $d_k =\frac{1}{r}(k\widetilde{d} -s_{k,max})$
and $s_k(E_0,F_0) = s_{k,max}.$ Choose inductively for any $\nu
=1,...,N_k$ an elementary transformation $E_{\nu}$ of $k$-type $I$ of
$E_{\nu -1}.$

In order to complete the proof of Theorem 3.2 it sufficies to show that
$$E_{N_k} \in U(r,d,k,s,i)$$ But $$s_k(E_{N_k}) =s_k(E_0) -N_kk =
s_{k,max} -N_kk=s$$ and $$
\begin{array}{lll}
s_i(E_{N_k})& \geq &s_i(E_0) -N_ki \  \ \ \ \  \ \ \mbox{ (by Proposition
1.6)}\\ &\geq&i(r-i)(g-1)-\frac{i}{k}(k(r-k)(g-1)-s+r-1)\\ &&\ \ \ \ \ \
\ \mbox{ (since $E_0$ is general and using (4).)}\\ &\geq
&i(k-i)(g-1)-\frac{i}{k}(r-2)\ \ \ \ \mbox{(since $s\geq 1$)}\\ &>&0 \ \ \
\ \ \ \  \ \ \mbox{(since $g\geq \frac{r+1}{2} $ by assumption)}
\end{array}
$$
\qed

\begin{rem} \label{3.3} \begin{em}
\begin{itemize}
\item[(a)] The assumption on the genus $g$ in 
Theorems 3.1 and 3.2 is imposed by the last line in the proof of Theorem
3.2.  The bound $g\geq \frac{r+1}{2}$ works for any $r$, for all $k, 1\leq
k\leq r-1$ simultaneously. If one fixes also $k$, the bound is slightly
better. In fact, if $k=1$ or $r-1,$ the proof shows that Theorems 3.1 and
3.2 are valid for any $g\geq 2.$ (For the proof note that in both cases
using duality one only has to check that $i(r-1-i)(g-1) - \frac{i}{r-1}
(r-2) >0$ for all $1 \leq i \leq r -2$. But this is valid for all $g \geq
2$.) For $2\leq k\leq r-2$ denote $k=\frac{r\pm n}{2}$ with $0\leq n\leq
r-4.$ Then the Theorems are valid for any $g\geq 3+2\frac{n-1}{r-n}.$
\item[(b)] There is a modification of the proof, for which the bound for
$g$ is also
slightly better. The duality can also be used to reduce the proof to the
case $k\geq \frac{r}{2}.$ Then one has also to show that $s_i(E_{N_k}) >0$
for $k<i<r.$ For this one has to choose the sequence of bundles
$E_0,E_1,...,E_{N_k}$ more carefully: Whenever possible one should use an
elementary transformation of $k$-type $I$ which is of $i$-type $II.$
\end{itemize}
\end{em} \end{rem}

\section{Stratification of ${\cal M}(r,d)$ according to the invariant
${{\se}_k}$} 

The function ${}{\se_k}: {\cal M}(r,d)\to\z$ defined by $E\mapsto
{}{\se_k}(E)$ is lower semicontinuous and this induces a stratification of
the moduli space ${\cal M}(r,d)$ into locally closed subvarieties $${\cal
M}(r,d,k,s):= \{ E\in {\cal M}(r,d) : s_k(E)=s \}$$ according to the value
$s$ of $s_k$. It is not clear (to us) whether ${\cal M}(r,d,k,s)$ is
irreducible or consists of several components.  Consider the natural map
$$ \phi :U(r,d,k,s) \ra {\cal M}(r,d,k,s) \subset {\cal M}(r,d).$$ As an
image of an irreducible variety $Im\phi $ is irreducible.  Let ${\cal
M}^0(r,d,k,s)$ denote the Zariski closure of $Im\phi $ in ${\cal
M}(r,d,k,s).$

\begin{lemma}\label{4.1} ${\cal M}^0(r,d,k,s)$ is an irreducible component
of ${\cal M}(r,d,k,s),$ if is nonempty.
\end{lemma}
\noindent{\bf Proof:} According to \cite{nr}, Proposition 2.4 there is a
finite \'etale
covering $\pi : \widetilde{M} \ra {\cal M}(r,d,k,s) \subset {\cal M}(r,d)$
and a vector bundle ${\cal E} $ on $C\times \widetilde{M}$ such that $id
\times \pi $ is just the classifing map. Let $Q_{\cal E}$ denote the Quot
scheme of ${\cal E}$ and $$0\ra {\cal F} \ra p^*{\cal E} \ra {\cal G} \ra
0  \eqno(5)$$ the universal exact sequence on $C\times \widetilde{M}
\times Q_{\cal E}$. Here $p:C\times \widetilde{M} \times Q_{\cal E} \ra
C\times \widetilde{M}$ denotes the projection map. Certainly there are
finitely many components of $Q_{\cal E}$, the union of which we denote by
$Q_{\cal E}^{r-k,d-d_k}, $ such that for all closed points $(e,x) \in
\widetilde{M} \times Q_{\cal E}^{r-k,d-d_k}$ the restriction ${\cal
F}_{|_{C\times \{(e,x)\}}}$ is a maximal subbundle of rank $k$ and degree
$d_k =\frac{1}{r}(s-kd) $ of $E={\cal E}_{|_{C\times \{(e,x)\} }}$ and
moreover every maximal subbundle occurs as a restriction of the exact
sequence $(5)$ to $C\times \{(e,x)\} $ for some $(e,x) \in \widetilde{M}
\times  Q_{\cal E}^{r-k,d-d_k}.$ As in section 2, the maximality condition
implies that ${\cal F}$ and ${\cal G}$ are vector bundles. Since
stableness is an open condition, it follows that the set of points $(e,x)
\in \widetilde{M} \times  Q_{\cal E}^{r-k,d-d_k}$ such that ${\cal
F}_{|_{C\times \{(e,x)\}}}$ and ${\cal G}_{|_{C\times \{(e,x)\}}}$ are
stable consists of whole components of $\widetilde{M} \times  Q_{\cal
E}^{r-k,d-d_k}.$ Hence the closure of the set of points $E\in {\cal
M}(r,d,k,s)$ which admits a stable maximal subbundle $F$ of rank $k$ such
that $E/F$ is also stable, consists also of whole components of ${\cal
M}(r,d,k,s)$. But since the set of such points $E$ of ${\cal M}(r,d,k,s)$
is just the irreducible set $Im\Phi$, this implies the assertion. \qed

It would be interesting to give an example for which ${\cal
M}^0(r,d,k,s)\not={\cal M}(r,d,k,s).$ For an example where $Im\phi \not=
{\cal M}(r,d,k,s)$ see Remark 4.5 below. The following theorem gives us
the dimension of ${\cal M}^0(r,d,k,s).$

\begin{thm} \label{4.2} 
 Let $r,d,k$, and $s$ be integers with $r\ge 2 $, $1\le k\le r-1$, $1\geq
s\geq k(r-k)(g-1 +(r-2)$ and $s\equiv kd\;\mod\; r.$ Suppose the genus of
$C$ is $g\ge \frac{r+1}{2}$.  Then ${\cal M}^0(r,d,k,s)$ is a non-empty
algebraic variety with $$\dim {\cal M}^0(r,d,k,
s)=\left\{\begin{array}{lcl} (r^2+k^2-rk)(g-1)+s+1 &  & s<k(r-k)(g-1)\\
&{\rm if}& \\ r^2(g-1)+1 & & s\ge k(r-k)(g-1)\end{array}\right. $$
\end{thm}

\noindent{\bf Proof:} If $s\ge k(r-k)(g-1)$, then $s =k(r-k)(g-1)+\epsilon_k$ 
where $\epsilon_k$ is the unique integer with $0\le\epsilon_k\leq r-1$ and
$s\equiv kd\;\mod\; r$. Then $Im\phi$ is open and dense in ${\cal
M}(r,d)$, which gives the assertion in this case (see [6]).

So we may assume that $s <k(r-k)(g-1)$. Consider the open set $U(r,d,k,s)$
in the variety $\p =\p(R^1_{p_{12_*}}(p^*_{02}{\cal F}^*_2\otimes
p^*_{01}{\cal F}_1)^*)$ of Section 3.  According to Theorem 3.2
$U(r,d,k,s)$ is non-empty, open and dense in $\p$. According to the
definitions of $U(r,d,k,s)$ and ${\se_k}$ the natural map
$$\phi:U(r,d,k,s)\longrightarrow {\cal M}^0(r,d,k,s)\subseteq {\cal
M}(r,d)$$ is dominant. We have to compute the dimension of ${\cal
M}(r,d,k,s)$.

Let $q\in U(r,d,k,s)$ be a general closed point. If $0\to F\to E\to G\to
0$ denotes the corresponding exact sequence, then $\phi(q)=E$ and
\cite{l1}, Lemma 4.2 says that $$\dim\phi^{-1}(E)\le h^0(F^*\otimes G).$$
On the other hand, $F$ and $G$ are general vector bundles in their
corresponding moduli spaces. Hence according to \cite{h2} Th\'eor\`em 4.6,
the vector bundle $F^*\otimes G$ is non-special, implying $$h^0(F^*\otimes
G)=0$$ since $\deg \ (F^*\otimes G)=s < k(r-k)(g-1)={\rm rk}(F^*\otimes
G)(g-1)$. Hence the generic fibre of $\phi$ is finite and thus
\begin{eqnarray*}
\dim {\cal M}^0(r,d,k,s)&=&\dim U(r,d,k,s)\\
&=&\dim \p(R^1p_{12*}(p^*_{02}{\cal F}^*_2\otimes p^*_{01}{\cal F}_1)^*)\\
&=&\dim {\cal M}(k,d_1)+\dim {\cal M}(r-k,d-d_1)+k(r-k)(g-1)+s-1
\end{eqnarray*}
with $d_1=\deg F={1\over r}(kd -s)$. Note that
$\p(R^1p_{12*}(p^*_{02}{\cal F}^*_2\otimes p^*_{01}{\cal F}_1)^*)$ is a
projective bundle of rank $k(r-k)(g-1)+s-1$ over a finite covering of
${\cal M}(k,d_1)\times {\cal M}(r-k, d-d_1)$ (see \cite{l1}, p. 455).
Hence
\begin{eqnarray*}
\dim {\cal M}^0(r,d,k,s)&=&k^2(g-1)+1+(r-k)^2(g-1)+1+k(r-k)(g-1)+s-1\\
&=&(r^2+k^2-rk)(g-1)+s +1. \hspace{5,2cm} \Box
\end{eqnarray*}

The fibres $\phi^{-1}(E)$ of the map $\phi:U(r,d,k,s)\to {\cal
M}^0(r,d,k,s)$ of the proof of Theorem 4.2 are exactly the sets of stable
maximal subbundles of $E$, whose quotient is also stabile. However,
maximal subbundles are not necessarily stable. Noting that a maximal
subbundle of a maximal subbundle of $E$ is also a subbundle of $E$ (and
similarly for quotient bundles) one easily shows

\bigskip
\begin{propn} \label{4.3} 
 Suppose $E\in {\cal M}(r,d)$ with ${\se_k}(E)=s$ for some $1\le k\le
r-1$.  Let $F$ be a maximal subbundle of rank $k$ of $E$. Then
\begin{itemize}
\item[(i)] ${\se_{\nu }}(F)\ge{1\over r}(k-\nu s)$ for all $1\le \nu \le k-1$
\item[(ii)] ${\se_{\nu }}(E/F)\ge{1\over r}((r-k)-(r-k-\nu )s)$ for all
$1\le \nu \le r-k-1$.
\end{itemize}
\end{propn}

So in particular for highs  ${\se_{\nu}}(F)$ or ${\se_{\nu }}(E/F)$ might
be very negative.  According to Lemma 2.1 there is a canonical
identification of the set ${}{M_k}(E)$ of maximal subbundles of $E$ with
the Quot-scheme $Q = Quot _E ^{r-k,{1\over r}((r-k)d+s)}$.  Hence there is
a universal subbundle ${\cal F}$ of $p^*E$ on the scheme $C\times
{}{M_k}(E)$, where $p:C\times {}{M_k}(E)\to C$ denotes the projection map.
Denote $${\widetilde{M_k}}(E):=\{ F\in {}{M_k}(E)|F\;\;{\rm and}\;\;
E/F\;{\rm stable}\}.$$ and by $\widehat{{M_k}}(E)$ the Zariski closure of
${\widetilde{M_k}}(E)$ in ${}{M_k}(E)$.  Applying the openness of
stability to the universal subbundle ${\cal F}$ and the universal quotient
bundle $p^*E/{\cal F}$ of $p^*E$ one deduces that $\widehat{M_k}(E)$
consists of whole irreducibility components of ${M_k}(E)$, namely exactly
of those components of ${{}{M_k}}(E)$ which contain a stable subbundle $F$
of $E$ such that $E/F$ is also stable. Moreover $\widetilde{M_k}(E)$ is
open in $\widehat{M_k}(E)$.

By definition we may canonically identify $${\widetilde{M_k}}(E)\to
\phi^{-1}(E).$$ This implies that
$$\dim{\widehat{M_k}}(E)=\dim\phi^{-1}(E).\eqno(6)$$ and we may use the
map $\phi$ to compute the dimension of $\widehat{{M_k}}(E).$ Let  $r,d,k$
and $s$ be integers as above and $g(C)\geq \frac{r+1}{2}$. According to
Theorem 4.2 the variety ${\cal M}^0(r,d,k,s)$ is non empty.

\bigskip
\begin{thm} \label{4.4}
 For a general vector bundle $E$ in ${\cal M}^0(r,d,k,s)$ we have $\dim
{\widehat{M_k}}(E)=\max(s -k(r-k)(g-1),0).$
\end{thm}
\bigskip

In particular, if $E$ is general in ${\cal M}(r,d)$ there is a unique
integer $\epsilon_k$ with $0\le\epsilon_k\le r-1$ and
$k(r-k)(g-1)+\epsilon_k\equiv kd\;\mod\; r$ and we have $$\dim
{\widehat{M_k}}(E)=\epsilon_k$$ If $s$ is not maximal value, i.e. $\se
<k(r-k)(g-1)$, then a general vector bundle in ${\cal M}^0(r,d,k,s)$
admits only finitely many stable maximal subbundles such that $E/F$ is
also stable.

\noindent{\bf Proof of Theorem 4.4:} The natural map 
$\phi :U(r,d,k,s)\to {\cal M}^0(r,d,k,s)$ is a dominant morphism of
algebraic varieties by the definition of $U(r,d,k,s)$ and ${\cal
M}^0(r,d,k,s)$. Let $q\in U(r,d,k,s)$ be a general point and $0\to F\to
E\to G\to 0$ be the corresponding exact sequence. According to \cite{l1}
Lemma 4.2 and Hirschowitz' Theorem (see \cite{h2} Th\'eor\`eme 4.6) we
have
\begin{eqnarray*}
\dim\phi^{-1}(E)&\le& h^0(F^*\otimes G)\\
&\le&\max(s-k(r-k)(g-1),0)
\end{eqnarray*}
So equation (6) implies the assertion for $s\le k(r-k)(g-1)$.  For
$s>k(r-k)(g-1)$ it suffices to show that the local dimension of
$\phi^{-1}(E)$ at $q$ is equal to $ s-k(r-k)(g-1)$. But Mori showed in
\cite{mo} that $$h^0(F^*\otimes G)-h^1(F^*\otimes G)\le \dim_qQ_{\cal E}$$
Again by \cite{h1} Th\'eor\`eme 4.6 the vector bundle $F^*\otimes G$ is
non-special implying $h^1(F^*\otimes G)=0$ and thus
$$\dim\widehat{M_k}(E)=\dim\phi^{-1}(E)=h^0(F^*\otimes G)=s-k(r-k)(g-1).
\eqno\Box   $$

\begin{rem}\label{4.5}\begin{em}
Take $r=3$, $d=1$, $k=2$, $s=2$ and $g\geq 2.$ In \cite{bgn} it was proved
that there are extensions $0\ra {\cal O}^2 \ra E\ra L\ra 0 $ of a line
bundle $L$ of degree $1$ by the trivial bundle ${\cal O}^2$ such that $E$
is stable.  Actually, such bundle $E$ is in ${\cal M}(3,1,2,2)$, since
$\mu (E) <1$ and hence $s_2(E) =2.$ However, for such bundles there is no
stable subbundles of degree $0$ and hence $E\notin Im\phi .$ Such bundles
$E$ are in the Brill-Noether locus ${\cal W}^k_{r,d}.$ An interesting
problem is to study the relation between the Brill-Noether loci ${\cal
W}^k_{r,d}$ and the ${\cal M}^0(r,d,k,s)$ varieties.
\end{em}\end{rem}

L. Brambila-Paz\\
Departamento de Matem\'aticas,\\
UAM - Iztapalapa\\
C.P. 09340.M\'exico, D.F.\\
E-mail: lebp@xanum.uam.mx

Herbert Lange\\
Mathematisches Institut der Universit\"at\\
Bismarckstra${\beta}$e 1{\small 1/2} D-8520\\
Erlangen, Germany.\\
E-mail: lange@mi.uni-erlangen.de

\end{document}